\documentclass[reqno]{amsart}%
\usepackage{amsmath}
\usepackage{amsfonts}
\usepackage{CJK}
\usepackage{amssymb}
\usepackage{graphicx}%
\usepackage{hyperref}
\usepackage{amscd}
\usepackage{mathrsfs}
\usepackage{bbm}

\newtheorem{thm}{Theorem}[section]

\theoremstyle{definition}

\theoremstyle{Conjecture}

\theoremstyle{remark}
\newtheorem{rem}{Remark}[section]
\theoremstyle{Example}

\newcommand{\be}{\begin{equation}}
\newcommand{\ee}{\end{equation}}
\newcommand{\bea}{\begin{eqnarray}}
\newcommand{\eea}{\end{eqnarray}}
\newcommand{\ben}{\begin{eqnarray*}}
\newcommand{\een}{\end{eqnarray*}}
\newcommand{\bet}{\begin{equation}
\begin{split}}
\newcommand{\eet}{\end{split}
\end{equation}}

\begin{document}
\title[equisingular approximations of quasi-plurisubharmonic functions]
{A remark on equisingular approximations of quasi-plurisubharmonic functions}

\author{Qi'an Guan}
\address{Qi'an Guan: School of Mathematical Sciences, and Beijing International Center for Mathematical Research,
Peking University, Beijing, 100871, China.}
\email{guanqian@amss.ac.cn}
\author{Zhenqian Li}
\address{Zhenqian Li: School of Mathematical Sciences, Peking University, Beijing, 100871, China.}
\email{lizhenqian@amss.ac.cn}

\thanks{The first author was partially supported by NSFC-11522101 and NSFC-11431013.}

\date{\today}

\begin{abstract}
In this note, we will present global equisingular approximations of quasi-plurisubharmonic functions with stable analytic pluripolar sets on compact complex manifolds.
\end{abstract}

\maketitle

\section{Introduction}\label{sec:introduction}

Let $X$ be a compact complex manifold of dimension $n$. A function $\varphi$ on X is said to be \emph{quasi-plurisubharmonic (quasi-psh)} if $\varphi$ is locally equal to the sum of a plurisubharmonic function and of a smooth function, that is, if its complex Hessian is bounded below by a $(1,1)$-form with continuous or locally bounded coefficients. The \emph{multiplier ideal sheaf} $\mathscr{I}(\varphi)$ associated to a quasi-psh function $\varphi$ is defined to be the sheaf of germs of holomorphic functions $f$ such that $|f|^2e^{-2\varphi}$ is locally integrable with respect to the Lebesgue measure, which is a coherent ideal sheaf.

We say that a quasi-psh function $\varphi$ has \emph{analytic singularities} if $\varphi$ can be written locally as
$$\varphi=\frac{c}{2}\log\big(\sum\limits_{k=1}^{m}|f_{k}|^2\big)+O(1),$$
where $c\in\mathbb{R}^+$ and $f_{k}$ are holomorphic functions.

In \cite{De_92}, Demailly proved a basic result on the approximation of psh functions by psh functions with analytic singularities via Bergman kernels. Here, we recall a global version as follows:

\begin{thm} \emph{(\cite{De_92}; see also \cite{De, De_Abel}).} \label{De_analytic}
Let $\varphi$ be a quasi-psh function on a compact Hermitian manifold $(X,\omega)$ of dimension $n$ such that $\frac{i}{\pi}\partial\bar\partial\varphi\geq\gamma$ for some continuous $(1,1)$-form $\gamma$. Then there is a sequence of quasi-psh functions $\varphi_m$ such that $\varphi_m$ has the same singularities as a logarithm of a sum of squares of holomorphic functions and a decreasing sequence $\varepsilon_m>0$ converging to 0 such that

$(a)$ $\varphi(x)\leq\varphi_m(x)\leq\sup\limits_{|\zeta-x|<r}\varphi(\zeta)+C\big(\frac{|\log
      r|}{m}+r+\varepsilon_m\big)$ with respect to coordinate open sets covering $X$. In particular, $\varphi_m$ converges to $\varphi$ and in $L^1(X)$ and

$(b)$ $\nu(\varphi,x)-\frac{n}{m}\leq\nu(\varphi_m,x)\leq\nu(\varphi,x)$ for every $x\in X$, where $\nu(\varphi,x)$ is
      the Lelong number of $\varphi$ at $x$;

$(c)$ $\frac{i}{\pi}\partial\bar\partial\varphi\geq\gamma-\varepsilon_m\omega$.
\end{thm}

In this note, by the strong openness property of multiplier ideal sheaves \cite{G-Z_open}, we will present the following equisingular approximations quasi-psh functions with stable analytic pluripolar set $A:=N(\mathscr{I}(\varphi))$, the zero set of multiplier ideal sheaf $\mathscr{I}(\varphi)$ associated to $\varphi$.

\begin{thm} \label{main1}
Let $\varphi$ be a quasi-psh function on a compact complex manifold $X$. Then, there exists a sequence $(\hat\psi_k)$ of quasi-psh functions with analytic singularities on $X$ such that

$(1)$ $(\hat\psi_k)$ converges to $\varphi$ almost everywhere;

$(2)$ $(``\mbox{equisingularity}")$ $\mathscr{I}(\hat\psi_k)=\mathscr{I}(\varphi)$ for every $k$;

$(3)$ $(``\mbox{stability}")$ the polar sets of $\hat\psi_k$ coincide with $A$ for all $k$.
\end{thm}

\begin{rem}
In \cite{Guan}, the first author proved that one cannot add the requirement to make $\hat\psi_k$ converge to $\varphi$ decreasingly.
\end{rem}

In \cite{DPS01}, the following ``equisingular" regularization process was established, which provides more ¡°equisingularity¡± than Theorem \ref{De_analytic} in the sense that the multiplier ideal sheaves are preserved.

\begin{thm} \emph{(\cite{DPS01}, see also \cite{De, De_Abel}).} \label{De_smooth}
Let $T=\alpha+i\partial\bar\partial\varphi$ be a closed $(1,1)$-current on a compact Hermitian manifold $(X,\omega)$, where $\alpha$ is a smooth closed $(1,1)$-form and $\varphi$ a quasi-psh function. Let $\gamma$ be a continuous real $(1,1)$-form such that $T\geq\gamma$. Then one can write $\varphi=\lim_{m\to+\infty}\varphi_m$, where

$(a)$ $\varphi_m$ is smooth in the complement $X\backslash Z_m$ of an analytic set $Z_m\subset X$;

$(b)$ $(\varphi_m)$ is a decreasing sequence, and $Z_m\subset Z_{m+1}$ for all $m$;

$(c)$ $\int_X(e^{-2\varphi}-e^{-2\varphi_m})dV_{\omega}$ is finite for every $m$ and converges to 0 as $m\to+\infty$;

$(d)$ $(``\mbox{equisingularity}")$ $\mathscr{I}(\varphi_m)=\mathscr{I}(\varphi)$ for all $m$;

$(e)$ $T_m=\alpha+i\partial\bar\partial\varphi_m$ satisfies $T_m\geq\gamma-\varepsilon_m\omega$, where
      $\varepsilon_m\to0$ as $m\to+\infty$.
\end{thm}

Note that in Theorem \ref{De_smooth} the pluripolar set $Z_m$ of $\varphi_m$ is larger and larger as $m$ tends to infinity. Combining the strong openness property of multiplier ideal sheaves and Theorem \ref{De_smooth}, we also present the following

\begin{thm} \label{main2}
Let $\varphi$ be a quasi-psh function on a compact complex manifold $X$. Then, there exists a sequence $(\hat\psi_k)$ of quasi-psh functions on $X$ such that

$(1)$ $(\hat\psi_k)$ is decreasing and convergent to $\varphi$;

$(2)$ $\hat\psi_k$ is smooth on $X\backslash A$ for all $k$;

$(3)$ $e^{-2\varphi}-e^{-2\hat\psi_k}$ is locally integrable for every $k$, which implies
      $\mathscr{I}(\hat\psi_k)=\mathscr{I}(\varphi)$.
\end{thm}

\section{Proof of main results}

To prove main results, the following strong openness property of multiplier ideal sheaves is necessary.

\begin{thm} \emph{(\cite{G-Z_open}).} \label{SOC}
Let $\varphi$ be a quasi-psh function on complex manifold $X$ and $K\subset X$ a compact subset. Then, there exists $\varepsilon_K>0$ such that for any $0\leq\varepsilon\leq\varepsilon_K$ we have
$$\mathscr{I}((1+\varepsilon)\varphi)|_K=\mathscr{I}(\varphi)|_K.$$
\end{thm}

We are now in a position to prove our main results.\\

\noindent{\textbf{\emph{Proof of Theorem} \ref{main1}.}}
Step 1. By Theorem \ref{De_analytic}, there exists a finite open covering $(U_k)$ of $X$ by coordinate balls and a sequence of quasi-psh functions $\varphi_m$ on $X$ such that on every $U_k$, we have
$$\varphi_{m}|_{U_k}=\frac{1}{2m}\log\sum_j|\sigma_{j,m}^{(k)}|^2+O(1),$$
where $(\sigma_{j,m}^{(k)})$ is an orthonormal basis of $\mathcal{H}_{U_k}(m\varphi)$, the Hilbert space of holomorphic functions $f$ on $U_k$ satisfying $\int_{U_k}|f|^2e^{-2m\varphi}d\lambda_n<\infty.$\\

Step 2 (see Proposition 2.1 in \cite{G-Z_Lelong1}). There exists a quasi-psh function $\tilde\varphi$ on $X$ such that

$(i)$ $e^{-2\varphi}-e^{-2\tilde\varphi}$ is locally integrable, which implies $\mathscr{I}(\varphi)=\mathscr{I}(\tilde\varphi)$ on $X$;

$(ii)$ $\tilde\varphi\in L_{loc}^{\infty}(X\backslash A)$.\\

Let $x_0\in X$ be any point. Without loss of generality, we assume $x_0\in U_k$. By Proposition 5.7 in \cite{De}, there exists $j_0>0$ and a neighborhood $V_0\subset\subset U_k$ of $x_0$ such that $\sigma_{1,1}^{(k)},...,\sigma_{j_0,1}^{(k)}$ generate $\mathscr{I}(\varphi)$ on $V_0$, and
$$\log\sum_j|\sigma_{j,1}^{(k)}|^2=\log\sum_{j=1}^{j_0}|\sigma_{j,1}^{(k)}|^2+O(1).$$ on $U_k$.
It follows from Theorem \ref{SOC} that there exists a real number $p_0>0$ satisfying $$\int_{V_0}|\sigma_{j,1}^{(k)}|^2e^{-2(1+\frac{1}{p})\varphi}d\lambda_n<\infty$$
for any $1\leq j\leq j_0$ and any $p\geq p_0$. Then, we obtain that on $V_0$,
\begin{equation*}
\begin{split}
&\int_{V_0}(e^{-2\varphi}-e^{-2\max\{\varphi,p\varphi_1\}})d\lambda_n=\int_{\{\varphi<p\varphi_1\}\cap V_0}
 e^{2\cdot\frac{1}{p}\varphi-2(1+\frac{1}{p})\varphi}d\lambda_n\\
\leq&\int_{\{\varphi<p\varphi_1\}\cap V_0}e^{2\varphi_1-2(1+\frac{1}{p})\varphi}d\lambda_n
 \leq\int_{V_0}e^{2\varphi_1-2(1+\frac{1}{p})\varphi}d\lambda_n\\
\leq&C\cdot\int_{V_0}\sum_{j=1}^{j_0}|\sigma_{j,m}^{(k)}|^2e^{-2(1+\frac{1}{p})\varphi}d\lambda_n<\infty
\end{split}
\end{equation*}
for some constant $C>0$.

Thus, we infer from the compactness of $X$ that for sufficiently large $p_0$, by taking $\tilde\varphi:=\max\{\varphi,p\varphi_1\}$ with $p\geq p_0$, we have $e^{-2\varphi}-e^{-2\tilde\varphi}$ is locally integrable and $\tilde\varphi\in L_{loc}^{\infty}(X\backslash A)$.\\

Step 3. As the discussion of Remark 3 in \cite{De_Abel}, by taking $\hat\varphi_m=(1+\frac{1}{m})\tilde\varphi_m$, we have $\mathscr{I}(\hat\varphi_m)=\mathscr{I}(\tilde\varphi)$ for large enough $m$, where $\tilde\varphi_m$ is the Bergman approximation sequence of $\tilde\varphi$. Moreover, $\hat\varphi_m$ is smooth outside the polar set $Z_m$ of $\hat\varphi_m$.

Since $\hat\varphi_m|_{X\backslash A}$ is smooth, we have $Z_m\subset A$. Then, it follows from $$\mathscr{I}(\hat\varphi_m)_x=\mathscr{I}(\tilde\varphi)_x\not=\mathcal{O}_{X,x}$$ for any $x\in A$ that $Z_m=A$ for large enough $m$.\\

Step 4. Equisingular approximation of $\varphi$.\\

Without loss of generality, we assume that both $\varphi$ and $\varphi_1$ are negative on $X$. Let $\psi_k:=\max\{\varphi,p_k\varphi_1\}$, where $p_0\leq p_k<p_{k+1}\to\infty\ (k\to\infty)$. Then $\psi_k$ converges to $\varphi$ and $\mathscr{I}(\psi_k)=\mathscr{I}(\varphi)$ for every $k$. By Remark 3 in \cite{De_Abel}, for every $k$, there exists a sequence of quasi-psh functions $\psi_{jk}$ with analytic singularities such that $\mathscr{I}(\psi_{jk})=\mathscr{I}(\varphi_k)$ for every $j$, and $\psi_{jk}$ is convergent to $\varphi_k$. It follows that for every $k$, there exists $N_k>0$ such that $m(E\{|\psi_{jk}-\psi_k|\geq\frac{1}{2^k}\})\leq\frac{1}{2^k}$, for any $j\geq N_k$.

Take a subsequence $\psi_{j_kk}$ of $\psi_{jk}$ with $N_k\leq j_k<j_{k+1}$. It follows from
$$|\psi_{j_kk}-\varphi|\leq|\psi_{j_kk}-\psi_k|+|\psi_k-\varphi|$$
that for every $\varepsilon>0$ and all $k$ with $\frac{1}{2^k}<\varepsilon$, we have
$$E\{|\psi_{j_kk}-\varphi|\geq\varepsilon\}\subset E\{|\psi_{j_kk}-\psi_k|\geq\frac{1}{2^k}\}\cup E\{|\psi_k-\varphi|\geq\varepsilon-\frac{1}{2^k}\},$$
which implies
$$m(E\{|\psi_{j_kk}-\varphi|\geq\varepsilon\})\leq\frac{1}{2^k}+\frac{\varepsilon}{2}<\varepsilon$$
for large enough $k$, i.e., $\psi_{j_kk}$ is convergent to $\varphi$ in Lebesgue measure. Let $\hat\psi_k$ be a subsequence of $\psi_{j_kk}$, which converges to $\varphi$ almost everywhere. Then, it follows from Step 3 that $\hat\psi_k$ satisfies (1), (2) and (3) as desired.
\hfill $\Box$
\\

\noindent{\textbf{\emph{Proof of Theorem} \ref{main2}.}}
Let $(\psi_k)$ be the equisingular approximation sequence of $\varphi$ as in Theorem \ref{De_smooth} and $\varphi_1, p_0$ as in Step 2 of Theorem \ref{main1}. In addition, we can assume all $\psi_k$ and $\varphi_1$ are negative. Then, we obtain that $e^{-2\psi_k}-e^{-2\varphi}$ and $e^{-2\varphi}-e^{-2\max\{\varphi,p\varphi_1\}}$ are locally integrable for any $p\geq p_0$, which implies that $e^{-2\psi_k}-e^{-2\max\{\varphi,p\varphi_1\}}$ is locally integrable. Since $$e^{-2\psi_k}-e^{-2\max\{\psi_k,\max\{\varphi,p\varphi_1\}\}}=
\max\{0,e^{-2\psi_k}-e^{-2\max\{\varphi,p\varphi_1\}}\}\in L_{loc}^1(X),$$
then $e^{-2\psi_k}-e^{-2\max\{\psi_k,p\varphi_1\}}$ is locally integrable.

Let
$$M_{\eta}(t_1,t_2):=\int_{\mathbb{R}^2}\max\{t_1+x_1,t_2+x_2\}
\prod\limits_{1\leq j\leq2}\eta_j^{-1}\theta(x_j/\eta_j)dx_1dx_2$$
be the regularized max function, where $\eta=(\eta_1,\eta_2)$ with $\eta_j>0$, and $\theta$ is a nonnegative smooth function on $\mathbb{R}$ with support in $[-1,1]$ such that $\int_{\mathbb{R}}\theta(x)dx=1$ and $\int_{\mathbb{R}}x\theta(x)dx=0$. By setting $\eta=(1,1)$, we have $M_{\eta}(\max\{\psi_k,p_0\varphi_1\},p\varphi_1)=\max\{\psi_k,p_0\varphi_1\}$ near $A$ for any $p>p_0$.
Then, it follows that
$$e^{-2\max\{\psi_k,p_0\varphi_1\}}-e^{-2M_{\eta}(\max\{\psi_k,p_0\varphi_1\},p\varphi_1)}$$
is locally integrable.

Note that
\begin{equation*}
\begin{split}
0\leq&e^{-2\psi_k}-e^{-2M_{\eta}(\psi_k,p\varphi_1)}\\
\leq&e^{-2\psi_k}-e^{-2M_{\eta}(\max\{\psi_k,p_0\varphi_1\},p\varphi_1)}\\
=&(e^{-2\psi_k}-e^{-2\max\{\psi_k,p_0\varphi_1\}})+
 (e^{-2\max\{\psi_k,p_0\varphi_1\}}-e^{-2M_{\eta}(\max\{\psi_k,p_0\varphi_1\},p\varphi_1)})\\
\end{split}
\end{equation*}
for any $p>p_0$. Thus, $e^{-2\psi_k}-e^{-2M_{\eta}(\psi_k,p\varphi_1)}$ is locally integrable for any $p>p_0$ and all $k$.

Let $\hat\psi_k:=M_{\eta}(\psi_k,p_k\varphi_1)$, where $\eta=(1,1)$ and $p_0<p_k<p_{k+1}\to\infty\ (k\to\infty)$.
Then, $\hat\psi_k$ decreasingly converges to $\varphi$ and smooth on $X\backslash A$. It follows from
$$e^{-2\varphi}-e^{-2\hat\psi_k}=(e^{-2\varphi}-e^{-2\psi_k})+(e^{-2\psi_k}-e^{-2M_{\eta}(\psi_k,p_k\varphi_1)})$$
that $e^{-2\varphi}-e^{-2\hat\psi_k}$ is locally integrable.
\hfill $\Box$

%


\begin{thebibliography}{123}
\bibitem{De_92} J.-P. Demailly, \emph{Regularization of closed positive currents and intersection theory}, J. Algebraic Geom. 1 (1992), 361--409.
\bibitem{De} J.-P. Demailly, \emph{Analytic Methods in Algebraic Geometry}, Higher Education Press, Beijing, 2010.
\bibitem{De_Abel} J.-P. Demailly, \emph{On the Cohomology of Pseudoeffective Line Bundles}, Complex Geometry and Dynamics, The Abel Symposium 2013, 51--99.
\bibitem{DPS01} J.-P. Demailly, T. Peternell, M. Schneider, \emph{Pseudo-effective line bundles on compact K\"ahler manifolds}, Internat. J. Math. 12 (2001), 689--741.
\bibitem{Guan} Q. A. Guan, \emph{Nonexistence of decreasing equisingular approximations with logarithmic poles}, J. Geom. Anal., published online on March 16, 2016; DOI: 10.1007/s12220-016-9702-2.
\bibitem{G-Z_open}  Q. A. Guan, X. Y. Zhou, \emph{A proof of Demailly's strong openness conjecture}, Ann. of Math. (2) 182 (2015), 605--616. See also arXiv:1311.3781.
\bibitem{G-Z_Lelong1} Q. A. Guan, X. Y. Zhou, \emph{Characterization of multiplier ideal sheaves with weights of Lelong number one}, Adv. Math. 285 (2015), 1688--1705.
\end{thebibliography}
\end{document}